\documentclass[12pt]{amsart}
\usepackage[margin=2.5cm]{geometry}
\usepackage{amssymb,color,url}
\normalsize
\def\bdm#1{\textcolor{black}{#1}}

\theoremstyle{definition}

\theoremstyle{remark}

\numberwithin{equation}{section}

\DeclareMathOperator{\PG}{\mathrm{PG}}

\DeclareMathOperator{\GF}{\mathrm{GF}}
\DeclareMathOperator{\PGammaL}{\mathrm{P}\Gamma{\mathrm{L}}}

\begin{document}


\title{There are 2834 spreads of lines in $\PG(3,8)$}


\author{Brendan D. McKay}
\address{Research School of Computer Science, 
Australian National University, Canberra ACT 0200, Australia}
\email{bdm@cs.anu.edu.au}


\author{Gordon F. Royle}
\address{Centre for the Mathematics of Symmetry and Computation, School of Mathematics and Statistics, University of Western Australia, 35 Stirling Highway, Nedlands WA 6009, Australia}
\email{gordon.royle@uwa.edu.au}



\begin{abstract}
In this note, we describe an exhaustive computer search for spreads of lines in $\PG(3,8)$ and determine that there are exactly
2834 inequivalent spreads under the group $\PGammaL(4,8)$. Therefore there are the same number of translation planes of 
order 64 with kernel containing $\GF(8)$, and we describe various properties of these planes.

Each portion of the search was performed at least twice with independently-written software, and the results checked for 
internal consistency by computation of the numbers of spreads rooted at a line not in the spread, 
thus enhancing confidence in the correctness of the search.
\end{abstract}


\maketitle

\section{Introduction}

A {\em line spread}, or just {\em spread} of the projective space $\PG(3,q)$ is a partition of the points of $\PG(3,q)$ into $q^2+1$ lines. It is well known that there is a {\em translation plane} of order $q^2$ associated with every line spread of $\PG(3,q)$, and that the planes associated to two spreads are isomorphic if and only if the two spreads are equivalent under the group $\PGammaL(4,q)$ of automorphisms of $\PG(3,q)$. If $q$ is {\em prime} then {\em every} translation plane of order $q^2$ arises from a spread of $\PG(3,q)$ and this property has been used to completely enumerate the translation planes of orders $q^2=25$ \cite{MR1149895} and $q^2=49$ \cite{MR1468940, MR1313377} by computing line spreads of $\PG(3,5)$ and $\PG(3,7)$ respectively.  In this note, we describe an enumeration of the line spreads of $\PG(3,8)$, and hence a certain class of translation planes of order $64$, and we describe various properties of these planes. 

In general, a {\em spread of $d$-spaces} is a partition of $\PG(2d-1,q)$ into disjoint copies of $\PG(d-1,q)$ and any spread of $d$-spaces yields a translation plane of order $q^d$. As $64 = 8^2 = 4^3 = 2^6$, this means that a translation plane of order $64$ may also arise from a partition of $\PG(5,4)$ into copies of $\PG(2,4)$ or a partition of $\PG(11,2)$ into copies of $\PG(5,2)$. Therefore our search is not a complete enumeration of the translation planes of order $64$, but merely those with kernel containing $\GF(8)$, which are likely to be a tiny fraction of the total number.

Let $\Gamma$ be the {\em line intersection graph} of $\PG(3,8)$; this is the graph with the $4745$ lines of $\PG(3,8)$ as its vertices, where collinear lines are adjacent. Then a {\em partial spread} of $\PG(3,8)$ is an independent set of $\Gamma$ and the maximum size independent sets of $\Gamma$ are precisely the spreads of $\PG(3,8)$. Thus the enumeration of spreads is equivalent to classifying the maximum independent sets in a particular graph, and as such can be tackled with the computational tools and techniques of graph theory.  

Of course, we 
are only interested in computing the spreads up to equivalence under $\PGammaL(4,8)$, the automorphism group of $\PG(3,8)$, and here care must be taken with the direct translation into a graph problem. This is because the automorphism group of 
$\Gamma$ is not equal to $\PGammaL(4,8)$, but rather it is twice as large. This happens because ${\rm Aut}(\Gamma)$ also
contains permutations of the lines determined by {\em dualities} of $\PG(3,8)$ (automorphisms that exchange points and hyperplanes). Thus two inequivalent spreads of $\PG(3,8)$ may in fact be equivalent under ${\rm Aut}(\Gamma)$. Thus, an equivalence class of maximum independent sets of $\Gamma$ can correspond to either $1$ or $2$ equivalence classes of spreads of $\PG(3,8)$.

Therefore we divide the process into the following two stages: The first stage is to use reasonably standard graph-theoretical techniques for computing one representative of each equivalence class of maximum independent sets of $\Gamma$ under ${\rm Aut}(\Gamma)$, and the second stage is to determine the representatives of equivalence classes of spreads under $\PGammaL(4,8)$.

The results, detailed in Section~\ref{results}, are that $\Gamma$ has exactly $1706$ equivalence classes
of independent sets, of which $578$ give a single equivalence class of spreads and $1128$ give two equivalence
classes of spreads, for a grand total of $578 + 2 \times 1128 = 2834$ equivalence classes of spreads of $\PG(3,8)$.


\section{Computational Details}

There are a number of standard techniques that take as input a graph $\Gamma$ and produce one representative of each equivalence class of independent sets of all sizes under ${\rm Aut}(\Gamma)$. In general, these proceed by augmenting smaller independent sets vertex-by-vertex and using either an explicit isomorphism check, or the implicit isomorphism checking used by orderly algorithms, to avoid constructing equivalent independent sets. 

In this situation however, we are only interested in the {\em maximum size} independents sets, and so a branch of the  search should be terminated as soon as it can be detected that the independent set currently under consideration has no completions of the desired size. However, a more serious problem is that an independent set can be constructed in this fashion by adding the vertices in any order, and there can be vast numbers of inequivalent {\em orderings} of the vertices, just as there are vast numbers of inequivalent partial spreads. 

To overcome these problems, we use the special structure of both the graph and the maximum independent sets that we are seeking.  As $\PG(3,8)$ has $585$ points and there are $73$ lines through each point, the edge set of $\Gamma$ consists of $585$ edge-disjoint cliques of size $73$, with each clique corresponding to all the lines through a particular point. Each pair of cliques meets in a unique vertex of $\Gamma$, namely the line connecting the corresponding pair of points.  As ${\rm Aut}(\Gamma)$ is transitive, we can freely choose any one line, say $\ell$, to be excluded from the spread. This determines nine cliques, one for each point of $\ell$, that each contribute a single vertex (other than $\ell$) to the spread. We use the term {\em starter}  or {\em starter based at $\ell$} to denote an independent set of size 9 consisting of one vertex from each of these cliques.

The computation of the inequivalent starters was independently performed twice, first by an orderly algorithm of the type described in \bdm{\cite{MR1606516, MR1614301}}, and secondly by a GAP program using the command {\tt SmallestImageSet} that, given a group and a $k$-set of points, computes the lexicographically least equivalent $k$-set.  This initial computation resulted in a collection of $1460$ starters. The correctness of this part of the computation can be verified theoretically. Given a line $\ell$ containing $9$ points $p_1$, $p_2$, $\ldots$, $p_9$, we can exactly count the total number of starters based at $\ell$: there are $72$ choices for the line through $p_1$, then $(72-8)$ choices for the line through $p_2$, then $72-(2\times 8)$ choices for the line through $p_3$, and so on. The number $72 - (k \times 8)$ of choices for the line through $p_{k+1}$
arises simply by counting the number of lines (other than $\ell$) joining $p_{k+1}$ to a point on one of the $k$ previously chosen lines. No line through $p_{k+1}$ other than $\ell$ can meet two of the previously chosen lines, or else all three lines would lie in a plane with $\ell$ and hence meet, and so these choices are distinct.  Thus counting the total number of pairs 
$(\ell, S)$ where $S$ is a starter based at $\ell$, we get
\[
4745 \times \prod_{j=0}^{j=8} (72 - 8 j) = \sum_S \frac{|G|}{|G_S|} t(S)
\]
where $S$ ranges over each of the $1460$ inequivalent starters, $G = {\rm Aut(\Gamma)}$, $G_S$ is the stabiliser of $S$ and $t(S)$ is the number of lines transversal to $S$ (that is, meeting every line of $S$).  As the computed values for the $1460$ starters satisfy this equation, we are confident that the collection of starters is correct and complete.

\medskip

The second stage of the computation is to process each starter individually, determining the spreads that arise from that particular starter.
\bdm{
The task can be described as solving a set of linear equations.
Define one variable $x_\ell$ for each of the 4745 lines, interpreted as 1 if $\ell$ is in
the spread and 0 otherwise.
Then, for each of the 585 points $p$, an equation $\sum_{p \in\ell} x_\ell=1$ states that exactly 
one line incident with $p$ is to be chosen.  Every solution for which all $x_\ell\in\{0,1\}$ is a spread.  
A starter $S$ fixes the values of some of the variables, either $x_\ell=1$ for those lines $\ell \in S$, and 
$x_m = 0$ for all the lines $m$ meeting $\ell$.
}

\bdm{
All solutions to the equation set for each starter were found twice.  One computation 
used an unpublished equation solver kindly provided by Petteri Kaski, while the other
was performed
using \textsc{minion}~\cite{minion}.  These took 8 years and 5 years of cpu time,
respectively, on a 
heterogeneous cluster of Linux workstations.  Fortunately, the results were identical. It is interesting to note that the
general purpose constraint satisfaction solver \textsc{minion} was faster than a highly optimised equation solver.
}

\bdm{
As the solutions were found, isomorphs under the action of ${\rm Aut}(\Gamma)$ were
removed, using \texttt{Traces}~\cite{NautyTraces}.  \texttt{Traces} is substantially
more efficient than \texttt{nauty} for this task, but since it was at the time experimental
software we took steps to verify it. Claims of isomorphism are safe since the isomorphism
is checked, while claims of non-isomorphism were verified by \texttt{nauty}  after
isomorphs were removed. All such verifications succeeded.
}

\bdm{
The nature of the search allowed a strong check at this point. Since our starters are
defined by some line $\ell$ which is not in the spread, the search should find at least
one member of each equivalence class of pairs $(\ell, S)$, where $\ell$ is a line and
$S$ is a spread that doesn't include~$\ell$.  We checked, for each pair $(\ell,S)$
that was discovered, that a pair equivalent to $(\ell',S)$ was also discovered for
each other line $\ell'$ not in~$S$.
This provides an additional check on the completeness of the starter set as well
as on the equation solving.
}

\medskip

The final stage of the computation is to determine whether the independent set represents one or two equivalence classes of spreads under $\PGammaL(4,8)$, which is a subgroup of index $2$ inside $G = {\rm Aut}(\Gamma)$. Let $g$ denote an arbitrary element of $G \ \backslash \PGammaL(4,8)$, and suppose that $S$ is a maximum independent set in $\Gamma$. If
\[
|G_S| = |\PGammaL(4,8)_S|
\]
then $S$ and $S^g$ are inequivalent under $\PGammaL(4,8)$, and so correspond to two non-isomorphic spreads, 
while if
\[
|G_S| = 2 |\PGammaL(4,8)_S|
\]
then $S$ and $S^g$ are equivalent under $\PGammaL(4,8)$, and yield a single spread.

\section{Results}\label{results}

There are 1706 pairwise inequivalent spreads of $\Gamma$ under the group ${\rm Aut}(\Gamma)$ which yield a total of $2834$ pairwise 
inequivalent spreads under $\PGammaL(4,8)$. The spectrum of automorphism group sizes is listed in Table~\ref{tab:autgroups}.

The data, namely the list of spread sets may be downloaded from the wRecall that a spread of $\PG(3,8)$ can be viewed as a collection of $2$-dimensional
subspaces of $\GF(8)^4$ that intersect only in the zero vector. Let
$W_\infty = \{(0,0,x,y) \mid x,y \in \GF(8)\}$ and
$W_0 = \{(x,y,0,0) \mid x,y \in \GF(8)\}$ be two such subspaces. 
The $2$-dimensional subspaces meeting these only in the zero vector all have the form
$W_A = \{(x, x A) \mid x \in \GF(8)^2\}$, where $A$ is a
non-singular $2 \times 2$ matrix and vectors are given as row-vectors. Under the action of $\PGammaL(4,8)$, any spread is equivalent to one
containing $W_\infty$ and $W_0$, which therefore has the form
\[
\{ W_\infty, W_0 \} \cup \{ W_{A_i} \mid 2 \leq i \leq 64 \}
\]
where each $A_i$ is a $2 \times 2$ matrix. The set of 64 matrices $\{A_1 = 0, A_2, A_3, \ldots, A_{64}\}$ 
is called a {\em spread set} and provides a compact description of any spread (up to equivalence). 
The condition that the subspaces $W_{A_i}$ and 
$W_{A_j}$ are disjoint can readily be seen to be equivalent to the
condition that $A_i - A_j$ is non-singular.  If $S$ is a spread with spread set $\{A_i\}$, then the set $\{A_i^T\}$ obtained by transposing all the matrices is another spread set, and so corresponds to another spread that we denote $S^T$. While $S$ and $S^T$ are equivalent under ${\rm Aut}(\Gamma)$ they may not be equivalent under $\PGammaL(4,8)$.

A complete list of the spread sets is available either on \url{arxiv.org} or from the first author's website at \url{http://cs.anu.edu.au/~bdm/data/geometries.html}. Each spread set occupies one line of 256 characters consisting of the 64 $2\times 2$ matrices written row-by-row. Each character is in the range $\{0,1,\ldots, 7\}$ and corresponds to an element of $\GF(8)$ where $0$ represents $0$, and $j > 0$ represents $x^{j-1}$ where
$x$ is a primitive element of $\GF(8)$ satisfying $x^3+x+1=0$. One file contains the $1706$ spread sets pairwise inequivalent under ${\rm Aut}(\Gamma)$ and a second contains the $2834$ spread sets pairwise inequivalent under $\PGammaL(4,8)$.

\begin{table}
\begin{tabular}{|r|rrr||r|rrr|} 
\hline
Order&S.p&Not s.p&Total&Order&S.p&Not s.p&Total\\
\hline
1 & 240 & 1872 & 2112 & 2 & 108 & 240 & 348 \\
3 & 35 & 44 & 79 & 4 & 19 & 18 & 37 \\
5 & 0 & 2 & 2 & 6 & 60 & 28 & 88 \\
8 & 13 & 14 & 27 & 9 & 3 & 6 & 9 \\
10 & 1 & 0 & 1 &12 & 14 & 4 & 18 \\
14 & 0 & 4 & 4 & 15 & 1 & 0 & 1 \\
16 & 6 & 0 & 6 & 18 & 15 & 8 & 23 \\
21 & 1 & 0 & 1  & 24 & 11 & 8 & 19 \\
 27 & 1 & 2 & 3 & 32 & 2 & 0 & 2 \\
36 & 8 & 0 & 8 & 42 & 1 & 2 & 3 \\
48 & 3 & 0 & 3 & 54 & 3 & 0 & 3 \\
72 & 7 & 0 & 7 & 96 & 2 & 0 & 2 \\
108 & 2 & 0 & 2 &120 & 1 & 0 & 1 \\
128 & 0 & 2 & 2 & 168 & 1 & 0 & 1 \\
189 & 1 & 0 & 1 & 192 & 2 & 0 & 2 \\
216 & 2 & 0 & 2 & 324 & 2 & 0 & 2 \\
360 & 1 & 0 & 1 & 384 & 4 & 0 & 4 \\
486 & 1 & 0 & 1  &1152 & 1 & 0 & 1 \\
1344 & 1 & 0 & 1 & 1512 & 2 & 0 & 2 \\
4032 & 0 & 2 & 2 & 27216 & 1 & 0 & 1 \\
87360 & 1 & 0 & 1  &14152320 & 1 & 0 & 1 \\
\hline
\end{tabular}
\caption{Spreads by automorphism group size}
\label{tab:autgroups}
\end{table}

%
%


\section{Rank}

The rank, over $\GF(p)$, of the incidence matrix of a projective plane of order $p^m$ is a useful invariant of the plane (usually studied in
conjunction with other properties of the $\GF(p)$-linear code generated by the incidence matrix of plane) called the {\em $p$-rank} of the plane. A famous conjecture of Hamada \cite{MR0332515} asserts that the Desarguesian plane $\PG(2,p^m)$, which has $p$-rank 
\[
\binom{p+1}{2}^m + 1,
\]
has the lowest rank among all projective
planes of the same order. When $p=2$, $m=6$ this gives a rank of $730$ and as shown in Table~\ref{tab:2rank}, this is the lowest rank among this collection of translation planes by a considerable margin. Hamada originally made a broader conjecture --- applicable to a larger class of designs --- that has been shown to be false in general, but restricted to projective planes, the conjecture remains firmly open and seems well supported by the (admittedly limited) evidence available.

\begin{table}
\begin{center}
\begin{tabular}{|rr|rr|rr|rr|rr|}
\hline
$2$-rank&No.&$2$-rank&No.&$2$-rank&No.&$2$-rank&No.&$2$-rank&No.\\
\hline
730 & 1 &
898 & 1 &
922 & 1 &
994 & 1 &
1006 & 1 \\
1030 & 1 &
1042 & 2 &
1048 & 2 &
1051 & 1 &
1057 & 1 \\
1063 & 5 &
1066 & 1 &
1072 & 1& 
1078 & 2&
1090 & 1 \\
1096 & 3 &
1099 & 1 &
1102 & 4 &
1108 & 1 &
1111 & 2 \\
1114 & 1 &
1117 & 4 &
1120 & 8 &
1123 & 7 &
1126 & 6 \\
1129 & 4 &
1132 & 20 &
1135 & 30 &
1138 & 2721 &&\\
\hline
\end{tabular}
\end{center}
\caption{The $2$-ranks of the $2834$ planes of order $64$ in this collection}
\label{tab:2rank}
\end{table}

It is interesting to note that more than 96\% of the planes have $2$-rank equal to $1138$, but it is not clear whether this has any particular significance.

\newpage

\bibliographystyle{acm2url}
\bibliography{../../gordonmaster.bib}

\end{document}